\documentclass[11pt]{article}
\usepackage{graphicx}
\usepackage{amssymb}
\usepackage{epstopdf}
\def\R{{\bf R}}
\def\H{{\bf H}}

\def\a{\alpha}

\def\e{\epsilon}
\def\l{\lambda}

\newtheorem{theorem}{Theorem}

\newtheorem{lemma}[theorem]{Lemma}

\def\eproof{\rule{3mm}{3mm}\newline}

\title{Orthospectra of Geodesic Laminations and Dilogarithm Identities on Moduli Space}
\author{Martin Bridgeman}

\begin{document}

\maketitle
\begin{abstract}
Given a measured lamination $\lambda$ on a finite area hyperbolic surface we consider a  natural measure $M_{\lambda}$ on the  real line obtained by taking the push-forward of the volume measure of the unit tangent bundle of the surface under an intersection function associated with the lamination. We show that the measure $M_{\lambda}$ gives  summation identities  for the Rogers dilogarithm function on the moduli space of a surface.
\end{abstract}

\section{Introduction}

Let $S$ be a closed hyperbolic surface and $\lambda$ a geodesic lamination on $S$. We let $\Omega$ be the volume measure on the unit tangent bundle $T_{1}(S)$.  We let $\alpha(v)$ be the longest geodesic arc  containing $v$ as a tangent vector and which does not intersect $\lambda$ transversely in its interior. Generically $\alpha(v)$ will be a geodesic arc with endpoints on $\lambda$.

We define the function $L:T_{1}(S) \rightarrow \overline{\R}$ by letting $L(v) = \mbox{Length}(\alpha(v))$. We note that $L(v)$ is measurable but can be infinite. We define measure $M_{\lambda}$ on the real line by $M_{\lambda} = L_{*}\Omega$. Then $M_{\lambda}$ is a measure describing the distribution of the lengths of $\alpha(v)$.

We cut $S$ along $\lambda$ to obtain a surface with boundary denoted $S_{\l}$. A {\em $\lambda-$cusp} of $S$ is an ideal vertex of a component of $S_{\l}$.  We let $N_{\lambda}$ be the number of $\lambda-$cusps of $S$. We denote by $\{\alpha_{i}\}$ the  
geodesic arcs in $S_{\l}$ which have endpoints perpendicular to $\partial S_{\l} \subseteq \lambda$ and denote the length of $\alpha_{i}$ by $l_{i}$. We note that if a component of $S_{\l}$ is an ideal  $k-$gon then there are a finite number of geodesics $\alpha_{i}$ in this component. Otherwise there are an infinite number. We call the set $\{l_{i}\}$ (with multiplicities) the {\em $\lambda-$orthospectrum.}  By doubling $S-\lambda$ we see that the $\lambda-$orthospectrum corresponds to a subset of the closed geodesics of a finite area surface and therefore is a countable set.
 
We prove the following length spectrum  identity
\begin{equation}
 \sum_{i} {\cal L} \left(\frac{1}{\cosh^2\frac{l_{i}}{2}}\right) = \frac{\pi^{2}}{12}(6|\chi(S)| -N_{\lambda})
\label{spec}
\end{equation}
where $\cal L$ is a Rogers dilogarithm function (described below).

\section{Dilogarithms and Polylogarithms}
The $k^{th}$ polylogarithm function $\mbox{Li}_k$ is defined by the Taylor series
$$\mbox{Li}_{k}(z) = \sum_{i=1}^{\infty} \frac{z^n}{n^k }$$
for $|z| < 1$ and by analytic continuation to $C$. In particular
$$\mbox{Li}_{0}(z) = \frac{1}{1-z} \qquad \qquad
\mbox{Li}_{1}(z) = -\log(1-z).$$
Also 
$$\mbox{Li}'_{k}(z) = \frac{\mbox{Li}_{k-1}(z)}{z} \qquad \mbox{giving} \qquad \mbox{Li}_{k}(z) = \int_{0}^{z} \frac{\mbox{Li}_{k-1}(z)}{z}\ dz.$$
Also the functions $\mbox{Li}_{k}$ are related to the Riemann $\zeta$ function by $\mbox{Li}_{k}(1) = \zeta(k)$.

The dilogarithm function is the function $\mbox{Li}_2(z)$ and is given by
$$\mbox{Li}_{2}(z) = -\int_{0}^{z} \frac{\log(1-z)}{z}\ dz.$$
Below is a brief description of some properties of the dilogarithm function. They can all be found in 1991 survey "Structural Properties of Polylogarithms" by L. Lewin (see \cite{Lew91}). From the power series representation, it is easy to see that the dilogarithm function satisfies the functional equation
$$\mbox{Li}_{2}(z) + \mbox{Li}_{2}(-z) = \frac{1}{2}\mbox{Li}_{2}(z^{2}).$$
Other functional relations of the dilogarithm can be best described by normalizing the dilogarithm function.
The (extended) Rogers ${\cal L}-$function (see \cite{Rog07}) is defined by
$${\cal L}(x) = \mbox{Li}_{2}(x) + \frac{1}{2}\log|x|\log (1-x)\qquad x \leq 1.$$ 
In terms of the Rogers ${\cal L}-$function, Euler's reflection relations  for the dilogarithm are 
\begin{eqnarray}
{\cal L}(x) + {\cal L}(1-x) = {\cal L}(1) = \frac{\pi^2}{6} \qquad 0 \leq x \leq 1 \nonumber \\
 {\cal L}(-x) + {\cal L}(-x^{-1}) = 2{\cal L}(-1) = -\frac{\pi^{2}}{6}  \qquad x >  0 
 \label{euler}
\end{eqnarray}
Also in terms of ${\cal L}$, Landen's identity is 
\begin{equation}
{\cal L}\left(\frac{-x}{1-x}\right) = -{\cal L}(x)  \qquad  0 < x < 1
  \label{landen}
  \end{equation}
  and Abel's functional equation  is
  \begin{equation}
  {\cal L}(x) + {\cal L}(y) = {\cal L}(xy) + {\cal L}\left(\frac{x(1-y)}{1-xy}\right) + {\cal L}\left(\frac{y(1-x)}{1-xy}\right).
  \label{abel}
  \end{equation}
  
Also a closed form for ${\cal L}(x)$ is known for certain values of $x$ including
$${\cal L}\left(1\right) = \frac{\pi^2}{6} \qquad {\cal L}\left(\frac{1}{2}\right) = \frac{\pi^2}{12} \qquad\mbox{(Euler)} \qquad {\cal L}\left(\phi^{-1}\right) = \frac{\pi^2}{10} \qquad {\cal L}\left(1-\phi^{-1}\right) = \frac{\pi^2}{15} \qquad\mbox{(Landen)} \qquad$$
where $\phi$ is the golden ratio. 

\section{Statement of Results}
The main result of the paper is the following;

{\bf Main Theorem}
{\em There exists a function $\rho:\R^{2} \rightarrow \R$ such that infinitesimally
$$dM_{\lambda} = \left(\frac{4N_{\lambda}x^{2}}{\sinh^{2}{x}} + \sum_{i}\rho(l_{i},x)\right)dx$$
where $N_{\lambda}$ is the number of $\lambda-$cusps of $S.$ Furthermore  the total mass of the measure $\rho(l,x)dx$ on the real line  is given by
$$F(l) =  \int_{0}^{\infty} \rho(l,x)\ dx  = 8{\cal L}\left({\frac{1}{\cosh^2\frac{l_{i}}{2}}}\right) $$
In particular  the measure $M_{\lambda}$ depends only on the $\lambda-$orthospectrum. 
}

\section{Length Spectrum Identity}
As $M_{\lambda} = L_{*}\Omega$, $M_{\lambda}$ has total mass equal to the volume of $T_{1}(S).$ Therefore $M_{\lambda}(\R) = \Omega(T_{1}(S)) = 4\pi^{2}|\chi(S)|.$ Summing up the masses of  measures in the Main Theorem we immediately obtain the following.

{\bf  Length Spectrum Identity Theorem}
{\em  Let $\lambda$ be a geodesic lamination on a finite area hyperbolic surface $S$. Then the  $\lambda-$orthospectrum satisfies the following 
$$ \sum_{i} {\cal L} \left(\frac{1}{\cosh^2\frac{l_{i}}{2}}\right) = \frac{\pi^{2}(6|\chi(S)| -N_{\lambda})}{12}
$$
or equivalently
$$ \sum_{i} {\cal L} \left(-\frac{1}{\sinh^2\frac{l_{i}}{2}}\right) = \frac{\pi^{2}(6\chi(S) +N_{\lambda})}{12}
$$
}

By Landen's identity (see equation \ref{landen}) we have
$${\cal L}\left(\frac{1}{\cosh^{2}(\frac{l}{2})}\right) = -{\cal L}\left(-\frac{1}{\sinh^{2}(\frac{l}{2})}\right).$$
Thus we can see that the second form of the Length Spectrum Identity corresponds to the first via Landen's identity.

\newpage

\section{Length Spectrum Identity on Moduli Space}
We note that if $S$ is a connected hyperbolic surface of finite area with non-empty geodesic boundary, 
 letting $\lambda = \partial S$ then the Length Spectrum Identity  gives a summation identity on the Moduli space $Mod(S)$ of $S$. In this case the Euler characteristic $\chi(S)$ can be a fraction and is defined such that $2\pi \chi(S)$ is the negative of the area of $S$. This relation is an infinite relation except in the case when $S$ is an ideal polygon. In this case we will show that these finite identities include  the classical dilogarithm identities described above.

\subsection{Classical Identities and the Moduli space of ideal polygons}
For $S$ an ideal n-gon, the Length Spectrum Identity  is a finite summation relation. We will show that the associated relations give an infinite list of finite relations including the classical identities stated in the previous section.

If $\{l_i\}$ is a $\lambda-$orthospectrum, we will define two parameterizations by letting
$$ a_i = -\frac{1}{\sinh^2\frac{l_{i}}{2}} \qquad \qquad b_i = \frac{1}{\cosh^2\frac{l_{i}}{2}}.$$
We now consider the Poincar\'e disk model and let $x_{i}, i = 1,\ldots, n$ be the vertices  in anticlockwise cyclic ordering around the circle. Let $s_{i}$ be the side $x_{i}x_{i+i}$. Let $l_{ij}$ be the length of the diagonal between $s_{i}$ and $s_{j}$ for $|i-j| \geq 2$.
We define  the cross-ratio  by
$$[z_{1},z_{2},z_{3},z_{4}] = \frac{(z_{1}-z_{2})(z_{4}-z_{3})}{(z_{1}-z_{3})(z_{4}-z_{2})}.$$
As  the cross ratio is invariant  under M\"obius transformations,  we map 
the quadruple $(x_{i},x_{i+1},x_{j},x_{j+1})$ to $(-1,1,e^{l_{ij}},-e^{l_{ij}})$.Then
$$[x_{i},x_{i+1},x_{j},x_{j+1}] = [-1,1,e^{l_{ij}},-e^{l_{ij}}] = -\frac{(-1-1)(-e^{l_{ij}}-e^{l_{ij}})}{(-1-e^{l_{ij}})(-e^{l_{ij}}-1)} =   \frac{4e^{l_{ij}}}{(e^{l_{ij}}+1)^{2}} = \frac{1}{\cosh^{2}(\frac{l_{ij}}{2})} $$
As $S$ has area $(n-2)\pi$ and $n$ cusps, $\chi(S) = (n-2)/2$ and $N_{\lambda} = n$. Thus the Length Spectrum identity  becomes
\begin{equation}
\sum_{i,j} {\cal L}([x_{i},x_{i+1},x_{j},x_{j+1}]) =  \frac{(n-3)\pi^{2}}{6}
\label{idealpolygon}
\end{equation}
where the sum is over all ordered pairs $i,j$ such that the sides $s_{i},s_{j}$ are disjoint (at infinity). 
In terms of dilogarithms we get
\begin{equation}
\sum_{i,j}  \mbox{Li}_{2}\left([x_{i},x_{i+1},x_{j},x_{j+1}]\right) = \frac{(n-3)\pi^{2}}{6} -\frac{1}{2} \sum_{i,j} \log\left(1-[x_{i},x_{i+1},x_{j},x_{j+1}]\right)\log\left([x_{i},x_{i+1},x_{j},x_{j+1}]\right)
\end{equation}

\subsection{Some  Cases}

{\bf Quadrilateral:} The ideal quadrilateral has 4 cusps and two ortholengths $l_{1},l_{2}$. By elementary hyperbolic geometry we have $\sinh(l_{1}/2).\sinh(l_{2}/2) = 1$. Therefore $a_1.a_2 = 1$ and letting $a = a_1$ the Length Spectrum identity  is equivalent to the  the classical reflection identity of Euler.
\begin{equation}
{\cal L}(a) + {\cal L}(a^{-1})  = -\frac{\pi^{2}}{6}.
\label{quadmod}
\end{equation}
Also we have
$$b_2 = \frac{1}{\cosh^{2}(l_{2}/2)} = \frac{1}{1+ \sinh^{2}(l_{2}/2)} = \frac{1}{1+ \frac{1}{\sinh^{2}(l_{1}/2)}} = \frac{\sinh^2(l_{1}/2)}{\cosh^2(l_1/2)} = 1-\frac{1}{\cosh^2(l_1/2)} = 1-b_1
$$
Thus letting $b = b_1$, the Length Spectrum identity is equivalent to the Euler reflection identity
\begin{equation}
{\cal L}(b) + {\cal L}(1-b)  = \frac{\pi^{2}}{6}.
\label{quadmod1}
\end{equation}

{\bf Pentagon and Abel's Identity:} If we choose a general ideal pentagon then there are $5$ diagonals and therefore $5$ parameters $a_{i}$. We send three of the vertices to $0, 1, \infty$ and  the other two to $u,v$ with $0 < u < v < 1$. Then the cross ratios in terms of $u,v$ are 

$$u, \qquad 1-v, \qquad \frac{v-u}{v} , \qquad \frac{v-u}{1-u}, \qquad \frac{u(1-v)}{v(1-u)} $$
Putting into the equation we obtain the following equation.  
\begin{eqnarray}
{\cal L}\left(u\right)+{\cal L}\left(1-v\right)+{\cal L}\left(\frac{v-u}{v}\right)+{\cal L}\left(\frac{v-u}{1-u}\right)+ {\cal L}\left(\frac{u(1-v)}{v(1-u)}\right)  = \frac{\pi^{2}}{3}.
 \end{eqnarray}
Letting $x = u/v, y = v$, then we get 
\begin{eqnarray}
{\cal L}\left(xy\right)+{\cal L}\left(1-y\right)+{\cal L}\left(1-x\right)+{\cal L}\left(\frac{y(1-x)}{1-xy}\right)+ {\cal L}\left(\frac{x(1-y)}{1-xy}\right)  =  \frac{\pi^{2}}{3}.
 \end{eqnarray}
Now by applying Euler's reflection identities for $x,y$, we obtain Abel's identity for the Rogers ${\cal L}-$function.
\begin{eqnarray}
{\cal L}\left(x\right)+{\cal L}\left(y\right) = {\cal L}\left(xy\right)+{\cal L}\left(\frac{y(1-x)}{1-xy}\right)+ {\cal L}\left(\frac{x(1-y)}{1-xy}\right) .
 \end{eqnarray}

{\bf General equation:} We obtain similar finite identities in the general ideal $n-$gon case. In general we note that equation \ref{idealpolygon}  will have $(n-3)$ independent variables and will be given by the summation of evaluating ${\cal L}$ on $\frac{n(n-3)}{2}$ rational functions in the $(n-3)$ variables.

\subsection{Regular Ideal n-gon relation}
We now consider the dilogarithm equation for the specific case of a regular ideal n-gon. In this case the cross ratios can be calculated  and the dilogarithm formulas for specific values of the dilogarithm function. 

We consider a regular ideal $n-$gon in with center $0$ in the Poincar\'e disk model  and vertices at $v_{k} = u^{k}, k = 0,\ldots,n-1$ for $u = e^{\frac{2\pi i}{n}}$. Then equation \ref{idealpolygon} can be thought of as an equation on the roots of the polynomial $z^n = 1$. We have

$$[v_{0},v_{1}, v_{r},v_{r+1}] = -\frac{(1-u)(u^{r+1}-u^{r})}{(1-u^{r})(u^{r+1}-u)} = \frac{u^{r}(u-1)^2}{u.(u^{r}-1)^{2}} =  \frac{\sin^{2}(\frac{\pi}{n})}{\sin^{2}(\frac{r\pi}{n})}$$

For $r < n/2$ there are exactly $n$ distinct perpendiculars between sides separated by $r$ sides and for $r = n/2$ there are $n/2$ such sides. To take care of the even and odd case simultaneously we let $e_{n}$ be $1$ if $n$ is even and $0$ if $n$ is odd. Therefore we have
\begin{equation}
\sum_{r=2}^{\lceil{n/2}\rceil - 1} n.{\cal L}\left(\frac{\sin^{2}(\frac{\pi}{n})}{\sin^{2}(\frac{r\pi}{n})}\right)
+e_{n} \frac{n}{2}.{\cal L}\left(\sin^{2}
\left(\frac{\pi}{n}\right)\right)  =  \frac{(n-3)\pi^{2}}{6}
\label{finitemod}
\end{equation}
where $\lceil{x}\rceil$ is the least integer greater than or equal to $x$.  Dividing by $n$ we get
$$
\sum_{r=2}^{\lceil{n/2}\rceil -1} {\cal L}\left(\frac{\sin^{2}(\frac{\pi}{n})}{\sin^{2}(\frac{r\pi}{n})}\right)
+ \frac{e_{n}}{2}.{\cal L}\left(\sin^{2}\frac{\pi}{n}\right) =   \frac{(n-3)\pi^{2}}{6n}
$$

{\bf Limiting case:} We let $n$ go to infinity and obtain the equation 
\begin{eqnarray*}
\lim_{n\rightarrow \infty} \left( \sum_{r=2}^{\lceil{n/2}\rceil -1} {\cal L}\left(\frac{\sin^{2}(\frac{\pi}{n})}{\sin^{2}(\frac{r\pi}{n})}\right)
+ \frac{e_{n}}{2}.{\cal L}\left(\sin^{2}\frac{\pi}{n}\right) \right)=  \lim_{n\rightarrow \infty}  \frac{(n-3)\pi^{2}}{6n} =  \frac{\pi^{2}}{6} 
\end{eqnarray*}
This gives a Rogers ${\cal L}-$function series relation due to Lewin (see p. 298 of \cite{Lew91})
\begin{eqnarray*}
 \sum_{r=2}^{\infty} {\cal L}\left(\frac{1}{r^{2}}\right) =   \frac{\pi^{2}}{6} 
\end{eqnarray*}

{\bf Regular ideal quadrilateral:} This case is trivial $a_1 = a_2 = -1$, $b_1 = b_2  = 1/2$  and  equations \ref{quadmod},\ref{quadmod1}  give the  classical evaluations
$${\cal L}(-1) =  -\frac{\pi^{2}}{12} \qquad \mbox{and} \qquad {\cal L}(\frac{1}{2}) =  \frac{\pi^{2}}{12}.$$

{\bf Regular ideal pentagon, Golden Mean:} 
For the regular ideal  pentagon, the orthospectrum consists of 5 geodesics each of the same length $l$. Using the formula above for $n=5, r =2$ we obtain that $l$ satisfies
$$\cosh^{2}\left(\frac{l}{2}\right) = \frac{2}{\sqrt{5}+3} = \phi^2 $$
where $\phi$ is the golden mean. Therefore as $\phi^2 = \phi+1$ 
$$\sinh^{2}\left(\frac{l}{2}\right) =\phi^2 -1 = \phi$$ and
we have $a = - \phi^{-1}$. Thus the Length Spectrum Identity gives  the classical relations of Landen
 $${\cal L}(-\phi^{-1}) = -\frac{\pi^{2}}{15} \qquad  \qquad {\cal L}(\phi^{-2}) = \frac{\pi^{2}}{15}.$$
Applying the quadrilateral relations \ref{quadmod}, \ref{quadmod1} we also get 
$${\cal L}(-\phi) = -\frac{\pi^2}{6} - {\cal L}(-\phi^{-1}) = -\frac{\pi^2}{10} \qquad.$$

{\bf Regular ideal Hexagon:}
For a regular ideal hexagon, there are 9 elements of the orthospectrum, with the $6$ being perpendicular to sides one apart and three being perpendicular to opposite sides. Putting $n=6$ into equation \ref{finitemod}
 above then gives
$$6{\cal L}(\frac{1}{3}) + 3{\cal L}(\frac{1}{4}) = \frac{\pi^{2}}{2}$$

Before we prove the main theorem, we first consider the geometry of ideal quadrilaterals in the hyperbolic plane. 

\section{Intersections with ideal quadrilaterals}
Given two disjoint geodesics $g_1,g_2$ with perpendicular distance $l$ between them, let $Q$ be the ideal quadrilateral with opposite sides $g_1,g_2$. Then we can map $Q$ by a M\'obius transformation to the ideal quadrilateral $Q_{a}$ in the upper half-plane with vertices $a,0,1,\infty \in \overline{\bf R}$ where $a < 0$. Similarly we  can map $Q$ to the ideal quadrilateral $Q_{b}$ in the upper half-plane with vertices $0,b,1,\infty \in \overline{\bf R}$ where $b >  0$. Using cross-ratios we have that 
\begin{equation}
a = -\frac{1}{\sinh^{2}\frac{l_{a}}{2}}\qquad b = \frac{1}{\cosh^{2}\frac{l_{a}}{2}}
\label{a(l)}
\end{equation}

The choice of normalization $Q_{a}, Q_{b}$  leads to the equivalent forms of the Length Spectrum Identity. We choose normalization $Q_{a}$ for our calculations.

If $x,y \in {\bf R}, x \neq y$, we let $g(x,y)$ be the geodesic in the upper half plane with end points $x,y$. Then for $(x,y) \in (a, 0)\times (1,\infty)$, the geodesic $g(x,y)$ intersects $Q_{a}$ in a definite length denoted $L(x,y)$.

\vspace{.2in}
\begin{figure}
\begin{center}
\includegraphics[width=4in]{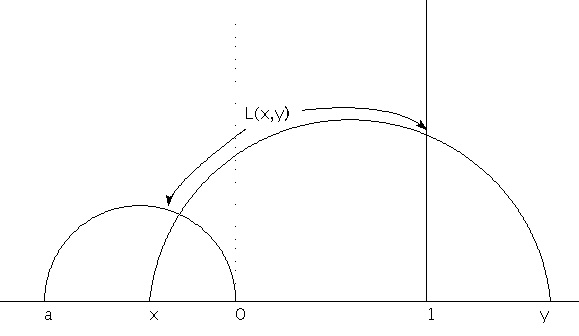}
\caption{Length intersection}
\end{center}
\end{figure}
\vspace{.2in}

\begin{lemma}
The map $L:(a,0) \times (1,\infty) \rightarrow {\bf R}$ is given by the formula
$$L(x,y) = \frac{1}{2} \ln\left( \frac{y(y-a)(x-1)}{x(x-a)(y-1)}\right) =\frac{1}{2} \ln{\frac{f(y)}{f(x)}} .$$
where $$f(x) = \frac{x(x-a)}{x-1}$$
\label{quad}
\end{lemma}

{\bf Proof:}
Let $T$ be the ideal triangle with vertices $0,1,\infty$. Let $l_{1}:(-\infty,0) \times (0,1) \rightarrow {\bf R}$ and  $l_{2}:(-\infty,0) \times (1,\infty) \rightarrow {\bf R}$ be given by letting $l_{1}(x,y)$ be the length of the intersection of $g(x,y)$ with $T$ and $l_{2}(x,y)$ be the length of the intersection of $g(x,y)$ with $T$.  By a previous paper (see \cite{BD07}) the functions $l_{i}$ are given by
$$l_{1}(x,y) =  \frac{1}{2}\ln\left(\frac{1-x}{1-y}\right) \qquad l_{2}(x,y) = \frac{1}{2}\ln\left(\frac{y(x-1)}{x(y-1)}\right)  $$

To calculate $L$, we split the quadrilateral $Q_{a}$ by the vertical line at $x = 0$ into two ideal triangles $T_{1}, T_{2}$ where $T_{1}$ has vertices $0,1, \infty$ and $T_{2}$ has vertices $a,0, \infty$.
Then $T_1 = T$ and  $f_{2}(z) = z/a$ sends $T_{2}$ to $T$. Therefore
$$L(x,y) = l_{2}(x,y) + l_{1}(y/a, x/a)$$

Therefore
$$L(x,y) = \frac{1}{2}\ln\left(\frac{y(x-1)}{x(y-1)}\right) +  \frac{1}{2}\ln\left(\frac{1-y/a}{1-x/a}\right) =  \frac{1}{2}\ln\left(\frac{y(x-1)(a-y)}{x(y-1)(a-x)}\right)
$$

\eproof
\vspace{.2in}
\begin{figure}
\begin{center}
\includegraphics[width=4in]{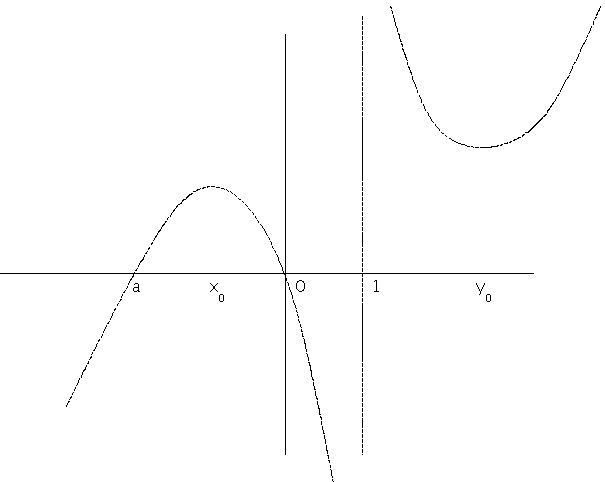}
\caption{Graph of function $f(x)$}
\end{center}
\end{figure}
\vspace{.2in}

We consider the rational  function $f(x)$ defined above. Differentiating we have
$$f'(x) = \frac{(2x-a)(x-1)-1.(x^{2}-ax)}{(x-1 )^{2}} = \frac{x^{2}-2x+a}{(x-1)^{2}}$$ 
Therefore  $f(x)$ has  two critical points $1 \pm\sqrt{1-a}$. We label the critical points $x_{0} = 1 - \sqrt{1-a}$ and $y_{0} = 1+\sqrt{1-a}$ and note that $x_{0}$ is a maximum and $y_{0}$ a minimum.

\section{Proof  of Summation Identity}
By definition
$$(L_{*}\Omega)(\phi) = \int_{T_{1}(S)}\phi(L(v))\ d\Omega.$$
Let $\alpha,\beta$ be two arcs in $S_{\l}$ with endpoints on $\partial S_{\l}$. Then we say $\alpha \sim \beta$ if they are homotopic relative to the boundary $\partial S_{\l}$.

We define the sets $A_{i} =\{v \in T_{1}(S) | \alpha(v) \sim \alpha_{i}\}$. Also for each $\lambda-$cusp $c$ we define $A_{c} = \{v \in T_{1}(S)| \alpha(v)  \sim c\}$ where $\alpha(v) \sim c$ if $\alpha(v)$ can be homotoped (rel boundary) out the cusp $c$.
Note that for $v \in A_{i}$ or $v \in A_{c}$, $L(v)$ is finite.
Finally we define the set $A_{\infty}$ to be all $v$ not in any $A_{i}$ or $A_{c}$.
By definition, the sets $A_{i}, A_{c}, A_{\infty}$ form a partition of $T_{1}(S)$. 
If we double $S_{\l}$ along its boundary, the geodesic arcs $\a_{i}$ correspond to a subset of the geodesics of the doubled surface.  Therefore as the length spectrum of the doubled surface is countable, so is the collection of arcs $\a_{i}$ in $S_{\l}$.
Also, by ergodicity of geodesic flow on $S$ (see \cite{Hop39}), the set $A_{\infty}$ is  a measure zero.

Therefore
$$(L_{*}\Omega)(\phi) = \sum_{i}\int_{A_{i}}\phi(L(v))\ d\Omega + \sum_{c} \int_{A_{c}}\phi(L(v))\ d\Omega.$$
 We let
$$a_{i} =  -\frac{1}{\sinh^{2}\frac{l_{i}}{2}}$$
Then setting $Q_{i} = Q_{a_{i}}$, we have that $Q_{i}$ is a quadrilateral with perpendicular of length $l_{i}$. 
We lift  $\alpha_{i}$ to the upper half plane so that it is the perpendicular of length $l_{i}$ in $Q_{i}$. 
 We lift each $\lambda-$cusp $c$ to the ideal vertex at infinity between the vertical geodesics $x = 0, x = 1$. Let $T$ be the ideal triangle with vertices $0,1,\infty \in \overline{\R}$.

If $v \in T_{1}(\H^{2})$ in the upper half plane, we define $g(v)$ to be the geodesic with tangent vector $v$. We also denote the endpoints of $g(v)$ by $(x(v), y(v))$.

We lift the set $A_{i}$ to the set   $A'_{i} \subseteq T_{1}(Q_{i})$. Then for $v \in A'_{i}$  the geodesic arc $\alpha'(v) = Q_{i}\cap g(v)$ is a lift of $\alpha(v)$. 
Similarly we lift $A_{c}$ to the set $A'_{c} \subseteq T_{1}(T)$.Then for $v \in A'_{c}$  the geodesic arc $\alpha'(v) = T \cap g(v)$ is a lift of $\alpha(v)$.
By abuse of notation we also let $\Omega$ be the volume measure on $T_{1}(\H^{2})$. We parameterize $T_{1}(\H^{2})$ by $(x,y,l) \in \overline{\R}\times\overline{\R}\times\R$ where $(x,y,l)$ corresponds to the vector $v$ such that $g(v)$ has ordered endpoints $(x,y)$ and $v$ has basepoint on $g(v)$ a distance $l$ from the highest point of $g(v)$ in the upper half-plane. Then the volume form $\Omega$ can be written as (see \cite{Nic89})
$$d\Omega = \frac{2dxdydl}{(x-y)^{2}}.$$
Therefore
$$\int_{A_{c}}\phi(L(v))\ d\Omega = \int_{A'_{c}}\frac{2.\phi(L(v))\ dxdydl}{(x-y)^{2}}$$

We note that $L(v)$ only depends on the endpoints and therefore we can write $L(v) = L(x,y)$. If $v \in A'_{c}$ then either $(x,y) \mbox{ or } (y,x) \in (-\infty, 0) \times (1,\infty)$ . Integrating over $l$ we have
$$\int_{A'_{c}} \frac{2.\phi(L(v)) \ dxdydl}{(x-y)^{2}} = \int_{-\infty}^{0} \int_{1}^{\infty}\frac{4.\phi(L(x,y))L(x,y)\ dxdy}{(x-y)^{2}}.$$
By our previous paper \cite{BD07}
$$ \int_{-\infty}^{0} \int_{1}^{\infty}\frac{4.\phi(L(x,y))L(x,y)\ dxdy}{(x-y)^{2}}. = \int_{0}^{\infty} \frac{4.\phi(L).L^{2}dL}{\sinh^{2}{L}}.$$
Thus as there are $N_{\lambda}$ $\lambda-$cusps we have
$$\sum_{c} \int_{A_{c}}\phi(L(v))\ d\Omega = N_{\lambda}. \int_{0}^{\infty} \frac{4.\phi(L).L^{2}dL}{\sinh^{2}{L}} =  M_{\infty}(\phi)$$
where $M_{\infty}$ is the measure with infinitesimal
$$dM_{\infty} = \frac{4N_{\lambda}x^{2}dx}{\sinh^{2}{x}}.$$
Similarly we have by lifting $A_{i}$ to $A'_{i}$ that
$$\int_{A_{i}}\phi(L(v))\ d\Omega = \int_{A'_{i}}\frac{2.\phi(L(v))\ dxdydl}{(x-y)^{2}}.$$
If $v \in A'_{i}$ then either $(x,y) \mbox{ or } (y,x) \in  (a_i,0) \times (1,\infty)$ . Integrating over $l$ we have
$$\int_{A'_{i}} \frac{2.\phi(L(v)) \ dxdydl}{(x-y)^{2}} = \int_{a_i}^{0} \int_{1}^{\infty}\frac{4.\phi(L(x,y))L(x,y)\ dxdy}{(x-y)^{2}}.$$
For $a < 0 $ we define $M_{a}(\phi)$ to be the righthandside of the above equation. Then 
$$M_{a}(\phi) = \int_{a}^{0} \int_{1}^{\infty}\frac{4.\phi(L(x,y))L(x,y)\ dxdy}{(x-y)^{2}}.$$
Then
$$M_{\lambda} = M_{\infty} + \sum_{i}M_{a_{i}}$$
As $M_{\lambda} = L_{*}\Omega$  it has total mass equal to the volume of $T_{1}(S)$ which is $4\pi^{2}|\chi(S)|$. Therefore
\begin{equation}
\Omega(T_{1}(S)) = 4\pi^{2}|\chi(S)| = M_{\lambda}(1) = M_{\infty}(1) + \sum_{i}M_{a_{i}}(1)
\label{volumesum}
\end{equation}
By an elementary calculation (see \cite{BD07})
$$\int_{0}^{\infty}\frac{x^{2}dx}{\sinh^{2}{x}} = \frac{\pi^{2}}{6}.$$
Therefore
$$M_{\infty}(1) = \int_{0}^{\infty} \frac{4N_{\lambda}x^{2}dx}{\sinh^{2}{x}} = 4N_{\lambda}\int_{0}^{\infty}\frac{x^{2}dx}{\sinh^{2}{x}} = 4N_{\lambda}.\frac{\pi^{2}}{6} = \frac{2N_{\lambda}\pi^{2}}{3}.$$
Using lemma \ref{quad} we substitute the formula for $L(x,y)$ to obtain
$$M_{a}(1) =  \int_{a}^{0} \int_{1}^{\infty}\frac{2.\log\left(\frac{y(y-a)(x-1)}{x(x-a)(y-1)}\right) dxdy}{(x-y)^{2}}.$$

Then by equation \ref{volumesum} above we obtain
$$4\pi^{2}|\chi(S)| = M_{\infty}(1) + \sum_{i}M_{a_{i}}(1) = \frac{2N_{\lambda}\pi^{2}}{3} + \sum_{i}F(l_{i})$$
giving the summation identity
\begin{equation}
 \sum_{i}F(l_{i}) = 4\pi^{2}|\chi(S)| - \frac{2N_{\lambda}\pi^{2}}{3} = \frac{2\pi^{2}}{3}(6|\chi(S)|-N_{\lambda})
 \label{functional}
 \end{equation}
 \eproof

\section{Integral Calculation}{\label{integral}}
In this section we find a formula for $F(l)$ by calculating an integral. We note that by the previous section, we already know that the function $F$ satisfies the functional equation \ref{functional}. We will make use of this to reduce $F$ to the form we wish independent of using any classical dilogarithm relations.

\begin{lemma}
For $ a < 0 $
$$\int_{a}^{0}\int_{1}^{\infty} \frac{\log\left|\frac{y(y-a)(x-1)}{x(x-a)(y-1)}\right|dxdy}{(x-y)^{2}} = -4{\cal L}(a)$$
\end{lemma}
{\bf Proof:} We let
$$G(a) =  \int_{a}^{0}\int_{1}^{\infty} \frac{\log\left|\frac{y(y-a)(x-1)}{x(x-a)(y-1)}\right|dxdy}{(x-y)^{2}}$$
Integrating by parts we get
$$\int \frac{\log\left|\frac{y(y-a)(x-1)}{x(x-a)(y-1)}\right|dx}{(x-y)^{2}} =
-\frac{\log\left|\frac{y(y-a)(x-1)}{x(x-a)(y-1)}\right|}{x-y}+\int \frac{1}{x-y}\left(\frac{1}{x-1} - \frac{1}{x}-\frac{1}{x-a}\right)dx.$$
Using
$$\int \frac{1}{(x-a)(x-b)} dx = \frac{1}{a-b}\left( \log|x-a| - \log|x-b|\right)$$ 
we get
\begin{eqnarray}
\int \frac{\log\left|\frac{y(y-a)(x-1)}{x(x-a)(y-1)}\right|dx}{(x-y)^{2}} =
-\frac{\log\left|\frac{y(y-a)(x-1)}{x(x-a)(y-1)}\right|}{x-y}+ \frac{1}{y-1}\left(\log|x-y| - \log|x-1|\right)+ \nonumber \\
 -\frac{1}{y}\left(\log|x-y|-\log|x| \right) -\frac{1}{y-a}\left(\log|x-y|-\log|x-a|\right) \nonumber \\
 = \frac{\log\left|\frac{y(y-a)(x-1)}{x(x-a)(y-1)}\right|}{y-x}- \frac{\log|x-1|}{y-1}
 +\frac{\log|x|}{y}  +\frac{\log|x-a|}{y-a}
 +\log|x-y|)\left(\frac{1}{y-1}- \frac{1}{y}- \frac{1}{y-a}\right) \nonumber
   \end{eqnarray}
We define
$$I(y)  =  \int_{a}^{0} \frac{\log\left|\frac{y(y-a)(x-1)}{x(x-a)(y-1)}\right|dx}{(x-y)^{2}} $$To evaluate the improper integral $I(y)$ we gather the divergent terms to find their limits. Therefore 
\begin{eqnarray}
I(y) =
 \lim_{x\rightarrow 0^{-}} \log|x|\left(\frac{1}{y}-\frac{1}{y-x}\right) - \lim_{x \rightarrow a^{+}} \log|x-a|\left(\frac{1}{y-a}-\frac{1}{y-x} \right) +\nonumber \\
\frac{\log\left|\frac{y(y-a)}{a(y-1)}\right|}{y} +\frac{\log|a|}{y-a}
 +\log|y|\left(\frac{1}{y-1}- \frac{1}{y}- \frac{1}{y-a}\right)\nonumber \\
 -\frac{\log\left|\frac{y(y-a)(a-1)}{a(y-1)}\right|}{y-a}
  +\frac{\log|a-1|}{y-1}-\frac{\log|a|}{y}
 -\log|a-y|\left(\frac{1}{y-1}- \frac{1}{y}- \frac{1}{y-a}\right) \nonumber
\end{eqnarray}
By elementary calculus, both limits are zero. As $y > 1$ and $a < 0$, when we gather the remaining terms by common denominators and get
\begin{eqnarray*}
I(y) = \frac{-2\log(-a)+2\log(y-a)-\log(y-1)}{y} + \frac{\log(1-a) + \log(y)-\log(y-a)}{y-1}+\\ +\frac{2\log(-a)-\log(1-a)-2\log(y)+\log(y-1)}{y-a}
\end{eqnarray*}
We now rewrite in the following form
\begin{equation}
I(y) = \left(\frac{\log(y)}{y-1}-\frac{\log(y-1)}{y}\right) +2\left(\frac{\log\left(\frac{y-a}{-a}\right)}{y}- \frac{\log\left(\frac{y}{-a}\right)}{y-a}\right) + \left(\frac{\log\left(\frac{y-1}{1-a}\right)}{y-a}-\frac{\log\left(\frac{y-a}{1-a}\right)}{y-1}\right).
\label{jform}
\end{equation}
Before we calculate the integral of $I(y)$ we note some properties of the dilogarithm.
As the dilogarithm function $\mbox{Li}_{2}$ satisfies 
$$\mbox{Li}_{2}(z) = -\int_{0}^{z} \frac{\log(1-t)}{t} dt$$
Then ${\cal L}$ has derivative
$${\cal L}'(x) = \frac{d}{dx}\left(\mbox{Li}_{2}(x) + \frac{1}{2}\log|x|\log(1-x)\right) = 
-\frac{\log(1-x)}{x}+\frac{1}{2}\left(\frac{\log(1-x)}{x} - \frac{\log|x|}{1-x}\right) = -\frac{1}{2}\left(\frac{\log(1-x)}{x}+\frac{\log|x|}{1-x} \right)$$
Now if $a < b$, then on the interval $x > b$, we have $(b-x)/(b-a) < 0$. We define
$$J(x,a,b) = 2{\cal L}\left(\frac{b-x}{b-a}\right).$$
Then differentiating $J$ we get
$$J'(x,a,b) = 2{\cal L}'\left(\frac{b-x}{b-a}\right).\frac{-1}{b-a} = \left( \frac{\log\left(\frac{x-a}{b-a}\right)}{b-x} +  \frac{\log\left(\frac{x-b}{b-a}\right)}{x-a}\right) = \left( \frac{\log\left(\frac{x-b}{b-a}\right)}{x-a} - \frac{\log\left(\frac{x-a}{b-a}\right)}{x-b} \right) $$
We set
$$J(y) = -J(y,0,1) -2J(y,a,0) + J(y,a,1).$$ Then  from equation \ref{jform} we have that
  $$J'(y) = -J'(y,0,1) -2J'(y,a,0) + J'(y,a,1) = I(y).$$
Therefore we have an antiderivative for $I$ and integrate to find $G$ to get
$$G(a) = \int_{1}^{\infty} I(y) dy = \left. J(y)\right|_{1}^{\infty} = \lim_{y \rightarrow \infty} J(y) - \lim_{y \rightarrow 1^{+}} J(y).$$ 
We let $ {\cal L}_{\infty}$ be the limit if ${\cal L}(x)$ as $x$ tends to $-\infty$. Therefore
$$\lim_{y \rightarrow 1^{+}} J(y) = -2{\cal L}(0) -4{\cal L}(a^{-1}) + 2{\cal L}(0) = -4{\cal L}(a^{-1}) \qquad\qquad \lim_{y \rightarrow \infty} J(y) = - 4{\cal L}_{\infty}.$$
Thus
$$G(a) =  -4{\cal L}_{\infty} + 4{\cal L}(a^{-1}) =-4({\cal L}_{\infty} -{\cal L}(a^{-1}))$$
It follows immediately from Euler's reflection identity that $G(a) = -4{\cal L}(a)$ but for completeness we derive it independently. From the formula we have $G(0) = -4 ({\cal L}_{\infty} -{\cal L}_{\infty}) = 0$. 
Also by equation \ref{functional}, $G(a)$ must satisfy a summation identity 
$$G(a) + G(a^{-1}) = G_{\infty}$$
where $G_{\infty}$ is the limit of $G$ as $a$ tends to $-\infty$.
Therefore 
$$G(a) = G_{\infty} - G(a^{-1}) = (G_{\infty} + 4{\cal L}_{\infty}) - 4{\cal L}(a)$$
But as $G(0) = 0$ we have
$$G(a) = -4{\cal L}(a)\qquad \mbox{and finally}\qquad F(l) =-8{\cal L}\left(-\frac{1}{\sinh^{2}(l/2)}\right).$$
\eproof

We note that by performing the integral over the  quadrilateral $Q_{b}$ where $b = 1/\cosh^{2}(l/2)$, the above can be repeated to show
$$F(l) =8{\cal L}\left(\frac{1}{\cosh^{2}(l/2)}\right).$$
Equivalently we note it also follows from Landen's identity.

\section{Volume interpretation of ${\cal L}$}
Let $g_{1}, g_{2}$ be disjoint geodesics in ${\bf H}^2$ with perpendicular distance $l$ and endpoints $x_1,y_1$ and $x_2,y_2$ respectively on ${\bf S}^1$. Given $v \in T_{1}(S)$ let $g_{v}$ be the associated oriented geodesic with tangent $v$. Then  we define the set 
$$C(g_1, g_2) = \left\{ v \in T_{1}(S) \left|\ g_{v} \cap g_1 \neq \emptyset,\   g_v \cap g_{2}\neq \emptyset \right\}\right.$$
Let $t = [x_1,y_1,x_2,y_2]$, then depending on the ordering of the points on the circle we have
$$ t = [-1,1,e^{l},-e^{l}] = \frac{1}{\cosh^{2}(l/2)} \qquad \mbox{or} \qquad t = [-1,1,-e^{l},e^{l}] = -\frac{1}{\sinh^{2}(l/2)}.$$
It follows  from the invariance of volume on $T_1(S)$, that the volume of $C(g_1, g_2)$ in $T_{1}(S)$ only depends on $t$. We therefore define $V(t) = \mbox{Volume}(S(g_1,g_2))$.
 
 Then it follows from the main theorem that
$${\cal L}(t) = \pm \frac{1}{8}V(t)$$
where the sign is given by the sign of $t$. Therefore we can interpret the Rogers ${\cal L}-$function  as a signed  volume function on $T_{1}(S)$ for the sets $G(g_1,g_2)$.

\section{Integral Formula for $\rho$}
We let 
$$L(x,y) =\frac{1}{2}\log\left( \frac{y(y-a)(x-1)}{x(x-a)(y-1)}\right) = \frac{1}{2}\log\left(\frac{f(y)}{f(x)}\right)\qquad\mbox{for}\qquad f(x) = \frac{x(x-a)}{x-1}.$$

Taking derivatives of the length function $L(x,y)$ we have
$$\frac{\partial L}{\partial x} = -\frac{f'(x)}{2f(x)} \qquad \frac{\partial L}{\partial y} = \frac{f'(y)}{2f(y)}.$$
By the previous section, the function $f$ has critical points $x_{0}, y_{0}$. Furthermore on $(a,0)$   the function $f(x)$ has global maximum at $x_{0}$ and  on $(1,\infty)$, $f$ has global minimum at $y_{0}$.
Therefore fixing $x$, the function $u:(1, \infty) \rightarrow {\bf R}$ given by $u(y) = L(x,y)$ is decreasing on $(1, y_{0})$ and increasing on $(y_{0}, \infty)$. Therefore we make the change of variable $t = L(x,y), x = x$. Finding inverses for $f$ we define the two function $g_{+}, g_{-}$ by
$$g_{\pm}(x) = \frac{(a+x) \pm\sqrt{(a+x)^{2} - 4x}}{2}.$$ 
Then solving $t = L(x,y)$ gives $f(y) = f(x)e^{2t}$. Therefore on $(1, y_{0})$ we have $y = g_{-}(f(x)e^{2t})$ and on $(y_{0},\infty)$ we have $y = g_{+}(f(x)e^{2t})$
Therefore
$$M_{a}(\phi) = \int_{a}^{0}\left( \int_{1}^{y_{0}} + \int_{y_{0}}^{\infty} \frac{4.\phi(L(x,y))L(x,y)\ dy}{(x-y)^{2}}\right)dx.$$
and
$$\int_{1}^{y_{0}} \frac{4.\phi(L(x,y))L(x,y)\ dy}{(x-y)^{2}} =  \int_{\infty}^{L(x,y_{0})}  \frac{4.\phi(t)t. g'_{-}(f(x)e^{2t})2f(x)e^{2t}dt}{(x-g_{-}(f(x)e^{2t}))^{2}}$$
$$\int_{y_{0}}^{\infty} \frac{4.\phi(L(x,y))L(x,y)\ dy}{(x-y)^{2}} =  \int_{L(x,y_{0})}^{\infty}  \frac{4.\phi(t)t. g'_{+}(f(x)e^{2t})2f(x)e^{2t}dt}{(x-g_{+}(f(x)e^{2t})^{2}}$$
Therefore combining  we have
$$M_{a}(\phi) = \int_{a}^{0} \int_{L(x,y_{0})}^{\infty}  8.\phi(t).t.e^{2t}.f(x) \left( \frac{g'_{+}(f(x)e^{2t})}{(x-g_{+}(f(x)e^{2t}))^{2}}-\frac{g'_{-}(f(x)e^{2t})}{(x-g_{-}(f(x)e^{2t}))^{2}}\right)dtdx$$

We switch the order of integration. The function $L(x, y_{0})$ is minimum at $x_{0}$ with minimum value $l = L(x_{0},y_{0})$ being the length of the perpendicular (see figure 2). Thus we integrate $t$ from $l$ to infinity. The integral in the $x$ direction is between the two $x$ solutions of $ t = L(x,y_{0})$ which are solutions to $f(x) = f(y_{0})e^{-2t}$. Thus we integrate $x$ from $g_{-}(f(y_{0})e^{-2t})$ to $g_{+}(f(y_{0})e^{-2t})$ giving
$$M_{a}(\phi) =  \int_{l_{a}}^{\infty}  8.\phi(t).t.e^{2t}dt  \left(\int_{g_{-}(f(y_{0})e^{-2t})}^{g_{+}(f(y_{0})e^{-2t})} \left(\frac{g'_{+}(f(x)e^{2t})}{(x-g_{+}(f(x)e^{2t}))^{2}}-\frac{g'_{-}(f(x)e^{2t})}{(x-g_{-}(f(x)e^{2t}))^{2}}\right)f(x)dx\right)$$
Therefore
$$M_{a}(\phi) =  \int_{0}^{\infty}  \phi(t).\rho(l,t) dt$$
where 
$$\rho(l,t) = 8t e^{2t}\chi_{[l, \infty)} . \left( \int_{g_{-}(f(y_{0})e^{-2t})}^{g_{+}(f(y_{0})e^{-2t})} \left(\frac{g'_{+}(f(x)e^{2t})}{(x-g_{+}(f(x)e^{2t}))^{2}}-\frac{g'_{-}(f(x)e^{2t})}{(x-g_{-}(f(x)e^{2t}))^{2}}\right)f(x)dx\right)$$
$$ \mbox{and } f(x) = \frac{x(x-a)}{x-1} \mbox{ where } a = - \frac{1}{\sinh^{2}{l/2}}$$

Therefore
$$(L_{*}\Omega)(\phi) = \int_{0}^{\infty} \phi(x) \rho(x)dx$$ where
$$\rho(x) = \frac{4N_{\lambda}x^{2}}{\sinh^{2}{x}} + \sum_{i} \rho(l_{i},x)$$
\eproof

\section{Asymptotic behavior}
In this section we study the asymptotic behavior of the function $\rho(l,t)$ for large $t$. 

For functions of a single variable, we write {\em $f(x) \simeq g(x)$ as $x$ tends to $x_{0}$} if 
$$\lim_{x \rightarrow x_{0}} \frac{f(x)}{g(x)} = 1.$$ 
Furthermore for functions of more than one variable, we write {\em $f(x,y) \simeq_{x} g(x,y)$ as $x$ tends to $x_{0}$} if 
$$\lim_{x \rightarrow x_{0}} \frac{f(x,y)}{g(x,y)} = 1.$$ 

\begin{theorem}
The measure $\rho(l,t)dx$ on the real line satisfies
$$   \lim_{t \rightarrow \infty}\frac{\rho(l,t)}{16t^{2}e^{-2t}} = r(l)$$
uniformly on compact subsets of  $(0,\infty)$ where $$r(l)  = \frac{-2a^{2}+5a-2}{a(1-a)} \qquad \mbox{ for }  \qquad a = -\frac{1}{\sinh^{2}\left(\frac{l}{2}\right)}$$
\end{theorem}

{\bf Proof:} We now show $\lim_{t \rightarrow \infty} \rho(l,t) = r(l)$ converges uniformly on compact subsets of $(0,\infty)$. Let $I \subseteq (0,\infty)$ be a compact interval. Now let $l \in I$. As before we let $a = -1/\sinh^{2}(l/2)$ and define $f(x) = x(x-a)/(x-1)$ with inverses $g_{\pm}$ and critical values $x_{0}, y_{0}$.
Let
$$G(t,x) = 8t e^{2t}\left(\frac{g'_{+}(f(x)e^{2t})}{(x-g_{+}(f(x)e^{2t}))^{2}}-\frac{g'_{-}(f(x)e^{2t})}{(x-g_{-}(f(x)e^{2t}))^{2}}\right)f(x)$$
Then for $t > l$ we have
$$\rho(l,t)  =   \int_{g_{-}(f(y_{0})e^{-2t})}^{g_{+}(f(y_{0})e^{-2t})} G(t,x)dx$$

For $C > 0$, we further define 
\begin{equation}
\rho(C,l,t) =   \int_{g_{-}(f(y_{0})Ce^{-2t})}^{g_{+}(f(y_{0})Ce^{-2t})} G(t,x)dx
\label{Capprox}
\end{equation}

On the interval $[a,0]$ $f$ has maximum at $x_{0}$. Therefore $\rho(C,l,t)$ is defined for all $t$ such that $f(y_{0})Ce^{-2t} < f(x_{0})$ or 
$$t  > K_{0}(C) =  \frac{1}{2}\ln{C} +   \frac{1}{2}\ln\left(\frac{f(y_{0})}{f(x_{0})}\right) =     l + \frac{1}{2}\ln{C}   $$ 
  
 Considering $g_{\pm}(x)$ for large $x$ we have
$$g_{\pm}(x) =  \frac{(a+x) \pm  
\sqrt{(a+x)^{2}-4x}}{2}   \simeq  \frac{(a+x)}{2}\left(1 \pm  
\left(1-\frac{2x}{(a+x)^{2}}\right)\right)$$
Therefore
$$g_{-}(x) \simeq \frac{(a+x)}{2}\left( 1 - 1  + \frac{2x}{(a+x)^{2}}\right) = \frac{x}{a+x} \simeq 1 - \frac{a}{x}$$
and
$$g_{+}(x) \simeq \frac{(a+x)}{2}\left( 1 + 1 - \frac{2x}{(a+x)^{2}}\right) = (a+x) - \frac{x}{a+x} \simeq (a-1) + x+\frac{a}{x}$$
Taking leading terms we have
\begin{equation}
 g_{-}(x) \simeq 1 \qquad \qquad g'_{-}(x) \simeq \frac{a}{x^2} 
\qquad\qquad\qquad g_{+}(x) \simeq x \qquad\qquad g'_{+}(x) \simeq 1
\label{gapprox}
\end{equation}

 We let $I_{C} = [g_{-}(f(y_{0})Ce^{-2t}), g_{+}(f(y_{0})Ce^{-2t})]$. Then for $x \in I_{C}$ we have
 $f(x)e^{2t} \geq C.f(y_{0})$. Therefore for $C$ sufficiently large we use the above approximations to approximate $G(t,x)$ on $I_{C}$.  We substitute the approximations  \ref{gapprox} into the formula for $G(t,x)$ to define
 $$G_{1}(t,x) = 8te^{2t}.\left(\frac{1}{(x-f(x)e^{2t})^{2}} - \frac{\frac{a}{(f(x)e^{2t})^{2}}}{(x-1)^{2}}\right)f(x)$$ Simplifying we have
 $$G_{1}(t,x) = 8te^{-2t}.\left( \frac{1}{(1- \frac{x}{f(x)e^{2t}})^{2}} - \frac{a}{(x-1)^{2}} \right)\frac{1}{f(x)}.$$
Noting that $f(x)e^{2t} > Cf(y_{0})$ on $I_{C}$, then for large $C$ the quantity  $\frac{x}{f(x)e^{2t}}$ is  small and we obtain the approximation 
$$G_{2}(t,x) = 8te^{-2t}.\left(1-\frac{a}{(x-1)^{2}}\right)\frac{1}{f(x)}.$$

Therefore given an $\e > 0$ we can find a $K_{1}(\e)$ such that
$$\frac{G(t,x)}{G_{2}(t,x)} \in [1-\e,1+\e] \qquad \mbox{ for all } C > K_{1}(\e), t > K_{0}(C), x \in I_{C}.$$
Therefore integrating
$$ \frac{1}{\rho(C,l,t)}  \left(  8te^{-2t} . \int_{g_{-}(f(y_{0})Ce^{-2t})}^{g_{+}(f(y_{0})Ce^{-2t})} 
\left(1- \frac{a}{(x-1)^{2}}\right)
 \frac{1}{f(x)} 
 dx 
 \right)  \in [1-\e,1+\e]$$
 for $C > K_{1}(\e)$ and $t > K_{0}(C) $.
We fix a $K > K_{1}(\e)$ and define 
$$\rho_{K}(l,t) = 8te^{-2t}.\left( \int_{g_{-}(f(y_{0})Ke^{-2t})}^{g_{+}(f(y_{0})Ke^{-2t})} \left(1-\frac{a}{(x-1)^{2}}\right) \frac{1}{f(x)} dx\right) $$
$$ = 8te^{-2t}.\left( \int_{g_{-}(f(y_{0})Ke^{-2t})}^{g_{+}(f(y_{0})Ke^{-2t})} \left(\frac{x-1}{x.(x-a)} -\frac{a}{x(x-a)(x-1)}\right)  dx\right) $$
Integrating we have
$$ \int \left(\frac{x-1}{x.(x-a)} -\frac{a}{x(x-a)(x-1)}\right)  dx =\left( \frac{1-a}{a}\ln{|x|} -\frac{a}{1-a}\ln{|x-1|} -\frac{a^2-3a+1}{a(1-a)} \ln{|x-a|)} \right)$$
Therefore
$$\rho_{K}(l,t) =  8te^{-2t}.\left( \frac{1-a}{a}\ln|x| -\frac{a}{1-a}\ln{|x-1|} - \frac{a^2-3a+1}{a(1-a)} \ln|x-a| \right)\left|_{g_{-}(f(y_{0})Ke^{-2t})}^{g_{+}(f(y_{0})Ke^{-2t})} \right.$$
For $x$ small we have
$$g_{\pm}(x) =  \frac{(a+x) \pm  
\sqrt{(a+x)^{2}-4x}}{2}   \simeq  \frac{(a+x)}{2}\left(1 \mp  
\left(1-\frac{2x}{(a+x)^{2}}\right)\right)$$
Therefore
$$g_{-}(x) \simeq (a+x) - \frac{x}{a+x} \simeq a - \frac{(1-a)x}{a}\qquad g_{+}(x) \simeq \frac{x}{a+x} \simeq \frac{x}{a} $$
Therefore
$$\rho_{K}(l,t) \simeq_{t} 8te^{-2t}. \left( \frac{1-a}{a}\ln\left|{\frac{Kf(y_{0})e^{-2t}}{a^{2}}}\right| -\frac{a}{1-a}\ln\left|\frac{1}{a-1}\right| - \frac{a^2-3a+1}{a(1-a)} \ln\left| \frac{a^{2}}{(1-a)f(y_{0})Ke^{-2t}}\right| \right)$$
Taking limits as we have
$$\rho_{K}(l,t) \simeq_{t} (16t^{2}e^{-2t}).\left(-\frac{1-a}{a} - \frac{a^2-3a+1}{a(1-a)} \right) =  (16t^{2}e^{-2t}).\frac{-2a^{2}+5a-2}{a(1-a)}$$
Therefore given $\epsilon > 0$ there exists $K_{1}(\epsilon)  > 0$ such that for any $C > K_{1}(\e)$  both
\begin{equation}
\liminf_{t \rightarrow \infty} \frac{\rho(C,l,t)}{16t^{2}e^{-2t}r(a)}
\ \  \mbox{and} \  \   \limsup_{t \rightarrow \infty} \frac{\rho(C,l,t)}{16t^{2}e^{-2t}r(a)} \mbox{ are in } [1-\e,1+\e].
\label{limsupliminf}
\end{equation}
where 
$$r(a) = \frac{-2a^{2}+5a-2}{a(1-a)}$$
We now define
$$\rho_{-}(C,l,t) =  \int_{g_{-}(f(y_{0})e^{-2t})}^{g_{-}(f(y_{0})Ce^{-2t})} G(t,x) dx \qquad\mbox{and}\qquad \rho_{+}(C,l,t) = \int_{g_{+}(f(y_{0})Ce^{-2t})}^{g_{+}(f(y_{0})e^{-2t})} G(t,x) dt.$$
Then by definition
$$\rho(l,t) = \rho(C,l,t) + \rho_{-}(C,l,t) + \rho_{+}(C,l,t).$$
We now bound the functions $\rho_{\pm}(C,l,t)$. Let $I^{-}_{C}, I^{+}_{C}$ be the given intervals.

On the interval $I$,  $g_{\pm}(f(x)e^{2t}) > 1$ and $x < 0$ so $(x-g_{\pm}(f(x)e^{2t}))^{2} > 1$. Also as $g'_{-}(f(x)e^{2t}) < 0$ we have
 $$|G(t,x)| = (8t.e^{2t}). \left(\frac{g'_{+}(f(x)e^{2t})}{(x-g_{+}(f(x)e^{2t}))^{2}} - \frac{g'_{-}(f(x)e^{2t})}{(x- g_{-}(f(x)e^{2t}))^{2}}\right)f(x).$$ 
 $$\leq 8t.e^{2t}.  \left(g'_{+}(f(x)e^{2t}) - g'_{-}(f(x)e^{2t})\right)f(x).$$ 
The derivative of $g_{\pm}(x)$ is given by
$$g'_{\pm}(x) = \frac{1}{2}\pm \frac{1}{2}\frac{x+a-2}{\sqrt{(a+x)^{2}-4x}}.$$
Therefore
$$g'_{+}(x) - g'_{-}(x) =  \frac{x+a-2}{\sqrt{(a+x)^{2}-4x}}.$$
As $f$ has critical values $f(x_{0})$ and $f(y_{0})$ we have that
$$g'_{+}(x) - g'_{-}(x) =  \frac{x+a-2}{\sqrt{(x-f(x_{0}))(x-f(y_{0}))}}.$$
We note that on $I^{\pm}_{C}$ we have $f(y_{0}) < f(x)e^{2t} < Cf(y_{0})$ then

$$g'_{+}( f(x)e^{2t}) - g'_{-}( f(x)e^{2t}) \leq  \frac{Cf(y_{0})+a-2}{\sqrt{(f(y_{0})-f(x_{0}))(f(x)e^{2t}-f(y_{0}))}}. $$
$$\leq  \left(\frac{Cf(y_{0})+a-2}{\sqrt{f(y_{0})-f(x_{0})}}\right)\frac{e^{-t}}{\sqrt{f(x)-f(y_{0})e^{-2t}}}$$
The function  $f(x) = x(x-a)/(x-1)$ has maximum at $x_{0}$ on $(a,0)$. Therefore  for $b < f(x_{0})$
$$f(x) - b = \frac{(x-g_{-}(b))(x-g_{+}(b))}{(x-1)}$$
As $x \in (a,0)$ we have
$$f(x)-b \geq (x-g_{-}(b))(g_{+}(b) - x)$$
Therefore
$$g'_{+}( f(x)e^{2t}) - g'_{-}( f(x)e^{2t})\leq  \left(\frac{Cf(y_{0})+a-2}{\sqrt{f(y_{0})-f(x_{0})}}\right)\frac{e^{-t}}{\sqrt{(x-g_{-}(f(y_{0})e^{-2t}))(g_{+}(f(y_{0})e^{-2t})-x)}}$$
Now restricting to $I^{+}_{C}$ we have $x > g_{+}(f(y_{0})Ce^{-2t})$. Therefore for $x \in I^{+}_{C}$,
$$g'_{+}( f(x)e^{2t}) - g'_{-}( f(x)e^{2t})\leq  \left(\frac{Cf(y_{0})+a-2}{\sqrt{(f(y_{0})-f(x_{0}))(g_{+}(f(y_{0})Ce^{-2t})-g_{-}(f(y_{0})e^{-2t}))}}\right)\frac{e^{-t}}{\sqrt{g_{+}(f(y_{0})e^{-2t})-x}}$$
Therefore we have
$$\rho_{+}(C,l,t) \leq \int_{I^{+}_{C}}|G(t,x)|dx \leq D(t)8te^{t}.\int_{I^{+}_{C}} \frac{f(x)}{\sqrt{g_{+}(f(y_{0})e^{-2t})-x}}dt$$
where $D(t)$ is the constant
$$D(t) =  \left(\frac{Cf(y_{0})+a-2}{\sqrt{(f(y_{0})-f(x_{0}))(g_{+}(f(y_{0})Ce^{-2t})-g_{-}(f(y_{0})e^{-2t}))}}\right)$$
As $f(x) = x(x-a)/(x-1)$ then, $0 < f(x) \leq ax$ on $(a,0)$ we have
$$\rho_{+}(C,l,t) \leq \int_{I^{+}_{C}}|G(t,x)|dx \leq D(t).8.a.t.e^{t}.\int_{I^{+}_{C}} \frac{x}{\sqrt{g_{+}(f(y_{0})e^{-2t})-x}}dx$$
By integration we have
$$\int_{a}^{b} \frac{x}{\sqrt{b-a}}dx = 
\frac{2}{3}(2b+a)\sqrt{b-a}$$
Therefore
$$\rho_{+}(C,l,t) \leq 16.D(t).a.t.e^{t}.\left(2g_{+}(f(y_{0})e^{-2t})+g_{+}(f(y_{0})Ce^{-2t})\right)\sqrt{g_{+}(f(y_{0})e^{-2t})-g_{+}(f(y_{0})Ce^{-2t})}.$$
Now for $t$ large we have
 $$\lim_{t\rightarrow \infty}D(t) = \left(\frac{Cf(y_{0})+a-2}{\sqrt{(f(y_{0})-f(x_{0})).|a|}}\right) = D.$$
We note for $x$ small $g_{+}(x) \simeq x/a$. Therefore 
$$\limsup_{t \rightarrow \infty} \left|\frac{\rho_{+}(C,l,t)}{t^{2}e^{-2t}} \right| \leq \limsup_{t \rightarrow \infty}\frac{16.D.a.t.e^{t}.\left(\frac{2f(y_{0})e^{-2t}+f(y_{0})Ce^{-2t}}{a}\right)
\sqrt{
\frac{f(y_{0})e^{-2t}- f(y_{0})Ce^{-2t}}
{a}}}{t^{2}e^{-2t}}.$$
$$\limsup_{t\rightarrow \infty} \left|\frac{\rho_{+}(C,l,t)}{t^{2}e^{-2t}} \right|  \leq \limsup_{t\rightarrow \infty}\frac{16.D.f(y_{0})^{3/2}(C+2)\sqrt{C- 1}}{t.\sqrt{-a}} = 0.$$
Thus
$$\lim_{t\rightarrow \infty} \frac{\rho_{+}(C,l,t)}{t^{2}e^{-2t}}  = 0.$$
Similarly for $\rho_{-}(C,l,t)$ we once again have that 
$$\lim_{t\rightarrow \infty} \frac{\rho_{-}(C,l,t)}{t^{2}e^{-2t}}  = 0.$$
Therefore given $\e > 0$ we can find  $K(\e)$ such that for $C > K(\e)$ by equations \ref{limsupliminf}
$$\limsup_{t \rightarrow \infty} \frac{\rho(l,t)}{16t^{2}e^{-2t}r(a)}  = \limsup_{t\rightarrow \infty} \left(\frac{\rho_{-}(C,l,t)}{16t^{2}e^{-2t}r(a)} + \frac{\rho(C,l,t)}{16t^{2}e^{-2t}r(a)} + \frac{\rho_{+}(C,l,t) }{16t^{2}e^{-2t}r(a)}\right) = \limsup_{t\rightarrow \infty}\frac{\rho(C,l,t)}{16t^{2}e^{-2t}r(a)} \in [1-\e,1+\e]$$
As $\epsilon$ is arbitrary we have
$$\limsup_{t \rightarrow \infty} \frac{\rho(l,t)}{16t^{2}e^{-2t}r(a)} = 1$$
Similarly
$$\liminf_{t \rightarrow \infty} \frac{\rho(l,t)}{16t^{2}e^{-2t}r(a)} = 1$$
\eproof


\end{document}